\documentclass{article}
\usepackage{shapepar,bm}
\usepackage[dvips,arrow,matrix,ps,line]{xy}
\usepackage{theorem,amsfonts,amssymb,amscd}
\usepackage{geometry}
\usepackage{makeidx}
\usepackage{stmaryrd}

\usepackage{graphicx}
\usepackage{graphics}

\usepackage{epstopdf}

\usepackage{amssymb,graphics,hyperref}

\hypersetup{
    colorlinks = true,
    linkcolor = blue,
    anchorcolor = red,
   citecolor = blue,
    filecolor = red,
    pagecolor = red,
    urlcolor = blue}

\newcommand{\mathsym}[1]{{}}

\makeindex

\def\wX{{[{\mathbf X}]}}

\newcommand{\llb}{\llbracket}
\newcommand{\rrb}{\rrbracket}
\def\bmtx{\begin{matrix}}
\def\emtx{\end{matrix}}
\def\NN{\mathbb N}

\def\bfz{{\mathbf z}}

\def\ovsig{\overline{\sigma}}

\def\V{V}

\def\wX{[{\bfX}]}

\def\bfu{{\mathbf u}}

\def\Qcal{\mathcal Q}

\def\d{\partial}
\def\Ucal{\mathcal U}
\def\Vcal{\mathcal V}

\def\bfx{{\mathbf x}}

\def\ZZ{\mathbb Z}

\def\CC{\mathbb C}

\def\QQ{\mathbb Q}

\def\PP{\mathbb P}

\def\cocoa{{\hbox{\rm C\kern-.13em o\kern-.07em C\kern-.13em o\kern-.15em A}}}

\def\Dcal{\mathcal D}

\def\bft{{\bf t}}

\def\bfu{{\bf u}}

\def\End{\mathrm{End}}

\def\blamb{{\bm \lambda}}
\def\bmu{{\bm \mu}}
\def\bnu{{\bm \nu}}
\def\Pcal{{\mathcal P}}

\def\w2M{\bigwedge^2M}

\def\w{\wedge }
\def\bw{\bigwedge }

\protect
\protect
\protect
\protect

\protect
\protect

\def\sra{\rightarrow}
\def\lra{\longrightarrow}

\def\proof{\noindent{\bf Proof.}\,\,}
\def\qed{{\hfill\vrule height4pt width4pt depth0pt}\medskip}
\def\be{\begin{equation}}
\def\ee{\end{equation}}
\def\bclm{\begin{claim}}
\def\eclm{\end{claim}}
\def\beqn{\begin{eqnarray}}
\def\eeqn{\end{eqnarray}}
\def\beqn*{\begin{eqnarray*}}
\def\eeqn*{\end{eqnarray*}}

\theoremstyle{change}
\theorembodyfont{\rmfamily}

\newtheorem{cor}{Corollary}[section]

\newtheorem{claim}{}[section]

\global
\global

\def\no@breaks#1{{\def\\{ \ignorespaces}#1}}    

  \makeatletter
\def\cleardoublepage{\clearpage\if@twoside \ifodd\c@page\else
\hbox{} \thispagestyle{empty}
\newpage
\if@twocolumn\hbox{}\newpage\fi\fi\fi} \makeatother

\usepackage{eso-pic,graphicx}
\makeatletter
\newcommand\BackgroundPicture[2]{%
  \setlength{\unitlength}{1pt}%
  default \put(0,\strip@pt\paperheight){%
  \parbox[t][\paperheight]{\paperwidth}{%
    \vfill
     \centering \includegraphics[angle=#2, width=15cm, height=15cm,  bb=0 0 150 150]{#1}
    \vfill
}}} %
\makeatother

\providecommand{\bysame}{\leavevmode\hbox to3em{\hrulefill}\thinspace}

\def\wX{{[\mathbf X]}}
\def\Exp{\mathrm{Exp}}

\title{Schur Polynomials and  Pl\"ucker degree of Schubert Varieties}
\author{Letterio Gatto}

\begin{document}

\maketitle
\begin{abstract}

\medskip
\noindent
The following is an informal report on the contributed talk given by the author during the INPANGA 2020[+1] meeting on Schubert Varieties.

The polynomial ring  $B$  in infinitely many indeterminates $(x_1,x_2,\ldots)$, with rational coefficients, has a vector space basis of Schur polynomials, parametrized by partitions. The goal of this note is to provide  an explanation of the following fact. If $\blamb$ is a partition of weight $d$,  then the partial derivative of order $d$ with respect to $x_1$ of the Schur polynomial $S_\blamb(\bfx)$  coincides with  the Pl\"ucker degree of the Schubert variety of dimension $d$ associated to $\blamb$, equal to the number of standard Young tableaux of shape $\blamb$. The generating function encoding all the degree of Schuberte varieties is determined and  some (known) corollaries are also discussed. 
\end{abstract}

\section{Introduction} 
\claim{} Let $G(r,n)$ be the complex Grassmann variety parametrizing $r$-dimensional vector subspaces of $\CC^n$, $\Qcal_r	\sra G(r,n)$ be its universal quotient bundle and $c_t(\Qcal_r)$ its Chern polynomial.  Following \cite[p.~271]{Ful}, let
$
F_\bullet(\blamb)\,\, :\,\, 0\subseteq F_1\subsetneq F_2\subsetneq \cdots\subsetneq F_r\subseteq \CC^n
$
be a flag of $r\geq 1$ subspaces of $\CC^n$ such that 
$
\dim F_i=i+\lambda_{r-i}.
$
Then $\blamb:=(\lambda_1\geq\cdots\geq \lambda_r)\in\Pcal_{r,n}$, a partition  whose Young diagram is contained in a $r\times (n-r)$ rectangle.
Let $\Omega^\blamb$ be the class of the closed (Schubert) irreducible variety of dimension $|\blamb|:=\lambda_1+\cdots+\lambda_r$ 
(\cite[Example 14.7.11]{Ful}).
$$
\Omega^\blamb(F_\bullet):=\Omega(F_1,\ldots,F_r):=\{\Lambda\in G(r,n)\,|\, 
\dim (\Lambda\cap F_i\geq i\},
$$ 
Its Pl\"ucker degree $f^\blamb$ (which coincides, by 
 \cite[Theorem 2.39]{smirnov}, with the number os standard Young tableaux of 
 shape $\blamb$),  does not depend on $n\geq r$.

\claim{} Let now $B:=\QQ[\bfx]$ be the polynomial ring in the infinitely many indeterminates $\bfx:=(x_1,x_2,\ldots)$. It  possesses a basis parametrized by  the set $\Pcal$ of all the partitions
\be
B:=\bigoplus_{\blamb\in \Pcal}\QQ\cdot S_\blamb(\bfx).
\ee 
If each indeterminate $x_i$ is given weight $i$, then  $S_\blamb(\bfx)$ is a homogeneous polynomial of weighted degree $|\blamb|$. Consider the vector subspace $\widetilde{B}_{r,n}:=\bigoplus_{\blamb\in\Pcal_{r,n}}\QQ\cdot S_\blamb(\bfx)$ of $B$ .
The map 
\be
\left\{\matrix{\pi_{r,n}&:&\widetilde{B}_{r,n}&\lra& H_*(G(r,n),\QQ)\cr\cr
&&S_\blamb(\bfx)&\longmapsto&\Omega^\blamb}\right.
\ee
is  a vector space isomorphism for trivial reasons.
The main result of this note is the following

\bclm{\bf Theorem.}\label{thm:thm13}{\em
\be
c_t(\Qcal_r)\cap \Omega^\blamb=\pi_{r,n}\Big(\exp\left(\sum_{i\geq 1}{t^i\over i}{\d\over \d x_i}\right)S_\blamb(\bfx)\Big)\label{eq1:mnth}
\ee
}
\eclm
In particular, equating  the coefficients of the  poer of $t$ of same degree:
\be
c_i(\Qcal_r)\cap \Omega^\blamb=S_i(\widetilde{\d})\cap \Omega^\blamb
\ee
where  $S_i(\widetilde{\d})$ iks an explicit polynomial expression in $\widetilde{\d}:=\displaystyle{\left({\d\over \d x_1}, {1\over 2}{\d\over \d x_2}, {1\over 3}{\d\over \d x_3},\ldots
\right)}
$
corresponding to  the coefficient of $t^i$ in the expansion of $\displaystyle{\exp\left(\sum_{i\geq 1}{t^i\over i}{\d\over \d x_i}\right)}$.
Theorem \ref{thm:thm13} will be shortly proven in Section \ref{sec:proof}, basing upon the notion of Schubert derivation on an exterior algebra as in \cite{SCHSD, gln, pluckercone}  alongwith its extension   to an infinite wedge power, as in \cite{SDIWP, BeGa}. 

\claim{} Theorem \ref{thm:thm13} has a number of corollaries, all collected in Section \ref{sec:sec3}. The most important is:

\smallskip
\noindent
{\bf Corollary \ref{cor:cor31}.}\label{thm:thm11} {\em For all $\blamb\in\Pcal$
\be
f^\blamb={\d^{d} S_\blamb(\bfx)\over \d x_1^{d}}.\label{eq10:fl}
\ee
}

\noindent
For example $f^{(2,2)}=\displaystyle{\d^4 S_{(2,2)}\over \d x_1^4}=2$, which is the Pl\"ucker degree of the Grassmannian $G(1,\PP^3)$ of lines in the three dimensional projective space.
We emphasize that we have not been able  to find any reference to (\ref{eq10:fl}), it looks new and is the main motivation of this note.

Recall that the classical way   to compute t$f^\blamb$ is to rely on a formula  due to Schubert, accounted for in \cite[Example 14.7.11]{Ful}. Because of its combinatorial interpretation in terms of Young tableaux,  it is also computed by the celebrated {\em hook length formula} (see e.g. 
 \cite[p.~53]{Fulyoung} or \cite[Theorem 4.33]{gillespie}) 
$$
f^\blamb:={|\blamb|!\over \prod_{x\in Y(\blamb)}h(x)},
$$
proved in \cite{thrall}, 
where $Y(\blamb)$ is the Young diagram of $\blamb$ and $h(x)$ is the hook 
length of the box $x\in Y(\blamb)$.

\medskip

\medskip
\noindent

\noindent
\smallskip
{\bf Corollary~\ref{cor:cor32}.} {\em 
\be
f^\blamb:={|\lambda|!\cdot \Delta_\blamb(\exp(t))}\label{eq:behzadf}
\ee }

\smallskip
\noindent
Formula (\ref{eq:behzadf}) has been first observed by O.~Behzad during the investigations which lead to her Ph.~D. Thesis \cite{BeThesis}. See also the forthcoming \cite{BeGa2}.

\medskip
\noindent
{\bf Corollary~\ref{cor:34}.} {\em Let $\blamb\in\Pcal$ and $Y(\blamb)$ its Young diagram. Then
\be
\prod_{x\in Y(\blamb)}h(x)={1\over \Delta_\blamb(\exp(t))}
\ee
where $h(x)$ denotes the {\em hook length} of the box $x$ in the Young diagram of $\blamb$.
}

\medskip
\noindent
{\bf Corollary \ref{cor:35}} {\em Let $\Pcal_r$ be the set of all partitions of length at most $r$. For all $\blamb\in\Pcal_r$, let $s_\blamb(\bfz_r)$ denote the Schur symmetric polynomial in the $r$ indeterminates $(z_1,\ldots,z_r)$, i.e.
$$
s_\blamb(\bfz_r)={\det(z_j^{\lambda_j-j+i})\over \Delta_0(\bfz_r)}
$$
Then
\be
\sum_{d\geq 0}{t^d\over d!}\sum_{\blamb\vdash n}f^\blamb s_\blamb(\bfz_r)=\exp(tp_1(\bfz))=\exp(t\cdot(z_1+\cdots+z_r))
\ee
In particular, for all $d\geq 0$
\be
(z_1+\cdots+z_r)^d=\sum_{\blamb\vdash d}f^\blamb\cdot s_\blamb(\bfz_r)\label{eqi:prewrds}
\ee
}
We additionally observe that evaluating the equality at $z_i=1$, formula (\ref{eqi:prewrds}) turns into
\be
r^d=\sum_{\blamb\vdash d}s_\blamb(1,\ldots,1)f^\blamb.\label{eq:comp}
\ee
Comparing (\ref{eq:comp}) with \cite[Formula (5), p.~52]{Fulyoung}, one deduces that  $s_\blamb(\underbrace{1,\ldots,1}_{r-times})$ is precisely the number $d_\blamb(r)$ of  standard Young tableaux  of shape $\blamb$,  whose entries are taken from the alphabet 
$\{1,2,\ldots,r\}.$

\claim{} Let $\displaystyle{\sum_{i\geq 0}}S_j(\tilde{\d})t^j=\exp\left(\sum_{i\geq 1}\displaystyle{t^i\over i}{\d\over \d x_i}\right)$. It is not difficult to see that
\be
\langle P(S_i(\bfx)),S_\blamb(\bfx)\rangle=P(S_i(\widetilde{\d}))S_\blamb(\bfx)
\ee
where by $P(S_i(\widetilde{\d}))$ is the evaluation of $P$ at $\displaystyle{x_i={1\over i}{\d\over \d x_i}}$.
In particular
$$
x_1^n=\sum_{\blamb\vdash n}<x_1^n,S_\blamb(\bfx)>S_\blamb(\bfx)=\sum_{\blamb\vdash n}{\d^n S_\blamb(\bfx)\over \d x_1^n}\cdot S_\blamb(\bfx),
$$
which, due to Corollary \ref{thm:thm11}, gives:
$$
x_1^n=\sum_{\blamb\vdash n}f_\blamb S_\blamb(\bfx)
$$
from which, taking the derivative with respect to $x_1$  of order $n$, and again by Corollary~\ref{cor:cor34} gives:
\be
n!=\sum_{\blamb\vdash n}(f^\blamb)^2\label{eq:square}
\ee
which is \cite[Formula (4), p.~50]{Fulyoung}.


\smallskip
\noindent

\noindent

In Section \ref{sec:sec2} we recall a few preliminaries. Then we will state and prove the main corollaries in Section \ref{sec:sec3}. In Section~\ref{sec:sec4} the notion of Schubert derivation (as in \cite{SCHSD}, \cite{HSDGA}, \cite{gln}) is reviewed. The short proof of Theorem \ref{thm:thm13} will conclude this short note.

\section{Preliminaries}\label{sec:sec2}
The content of this section is very well known and easily available in many 
common textbooks and its only  purpose is to introduce the notation adopted in the 
sequel.
\claim{\bf Partitions.} Let $\Pcal$ be the set of all partitions, namely the 
monoid of all non-increeasing sequences $\blamb:=(\lambda_1\geq \lambda_2\geq
\cdots)$ of non-negative integers with finite support (all terms zero but 
finitely many). The non zero terms of $\blamb$ are  called {\em parts}, the 
number $\ell(\blamb)$ of parts is called {\em length}. Let $\Pcal_r:=\Pcal
\cap \NN^r$: it is the  set of all partitions with at most $r$-parts and  $
\Pcal_\infty=\Pcal$.
The Young diagram of a partition is the left 
justified array of $r$-rows such that the $i$th row has $\lambda_i$ boxes. 
For all $\leq r\leq n$ we denote by $\Pcal_{r,n}$ the set of partitions whose 
Young diagram is contained in a $r\times (n-r)$ rectangle. Then $\Pcal_{r,
\infty}=\Pcal$ and $\Pcal_\infty=\Pcal$.
 If $\blamb\in 
\Pcal_{r,n}$, we denote by $\blamb^c$ the partition whose 
Young diagram is the complement of the Young diagram of $\blamb$ in the $r
\times (n-r)$ rectangle. For example the complement of the partition $(3,3,2,
1)$ in the $4\times 3$ rectangle is $(2,1)$. Its complement in  the $5\times 
4$ 
rectangle is $(4,3,2,1,1)$.
\claim{\bf Schur Determinants.}  Let $A$ be any commutative algebra. To each pair  
$$
\Big(f(t)=\sum_{n\geq 0}f_nt^n,\blamb\Big)\in A\llb t^{-1},t\rrb\times \Pcal_r
$$ 
one attachesthe {\em Schur determinant}:
\be
\Delta_\blamb(f(t))=\det(f_{\lambda_j-j+i})_{1\leq i,j\leq r}\in A\label{eq:scud}
\ee
If $f(t)\in A[[t]]$, one think of it as a formal Laurent series with $f_j=0$ for $j<0$.

%
Putting  $S_\blamb(\bfx):=\det(S_{\lambda_j-j+i})$, it is well known that
$$
B:=\bigoplus_{\blamb\in\Pcal}\QQ\cdot S_\blamb(\bfx)
$$


\claim{} \label{notpier} For all $(i,\blamb)\in\NN\times \Pcal_r$,  define
$$
PF_i(\blamb):=\{\bmu\in \Pcal_r\,\, |\,\, |\bmu|=|\blamb|+i\quad \mathrm{and}\quad \mu_1\geq\lambda_\geq\cdots\geq\mu_r\geq\lambda_r\}
$$
and
$$
PF_{-i}(\blamb):=\{\bmu\in \Pcal_r\,\, |\,\, |\bmu|=|\blamb|-i\quad \mathrm{and}\quad \lambda_1\geq \mu_1\geq\lambda_2\geq\mu_2\geq \cdots\geq\lambda_r\geq\mu_r\}
$$

\claim{\bf Proposition.} \label{prop:prop24}{\em Pieri's rule for Schur $S$-function holds: 
\be
S_i(\bfx)\cdot S_\blamb(\bfx)=\sum_{\bmu\in PF_i(\blamb)} S_\bmu(\bfx)\label{eq:PieriS}
\ee
}
\claim{\bf Proposition.}\label{prop:dualp} {\em The ``dual'' Pieri's rule for $\Omega^\blamb$ holds: 
\be
c_i(\Qcal_r)\cap \Omega^\blamb=\sum_{\bmu\in PF_{-i}(\blamb)} \Omega_\bmu(\bfx)\label{eq:PieriS}
\ee
}
Proposition \ref{prop:dualp} is basically another phrasing of \cite[Example 14.7.1]{Ful}
It is very well known that
$H_*(G(r,n),\QQ):=\bigoplus_{\blamb\in\Pcal_{r,n}}\QQ\cdot \Omega_\blamb$. Moreover,  Giambelli's formula says that
$$
\Omega^\lambda=\Delta_{\blamb^c}(c_t(\Qcal_r))\cap [G(r,n)]:=(c_{\lambda_j-j+i}(\Qcal_r))_{1\leq i,j\leq r},
$$
where $\blamb^c$ denotes the complement of $\blamb$ in the $r\times (n-r)$ rectangle.

\section{Corollaries to Theorem \ref{thm:thm13}}\label{sec:sec3}
In this section we assume Theorem \ref{thm:thm13}, to find a few corollaries. As it is customary to do we set $\sigma_i:=c_i(\Qcal_r)$.
\begin{cor}\label{cor:cor31} {Let $\blamb\in\Pcal_{r,n}$.
Then
\be
\sigma_1\cap \Omega^\blamb=\pi_{r,n}\left({\d S_\blamb\over \d x_1}\right)\label{eq:s1}
\ee
In particular, if $d=|\blamb|$:
\be
f^\blamb:=\sigma_1^d\cap \Omega^\blamb={\d^d\over\d x_1^{d}}S_\blamb(\bfx)\label{eq2:degr}
\ee
}
\end{cor}
\proof
Formula (\ref{eq:s1}) comes from equating the coefficients of the linear terms on both sides of~\ref{eq1:mnth}).
Iterating it $d$ times  one obtains (\ref{eq2:degr}), where the projection $\pi_{r,n}$ can be  omitted ($\pi_{r,n}$ is the identity on constants).\qed

Formula (\ref{eq:s1}) generalizes to
$$
\sigma_i\cap \Omega^\blamb=\pi_{r,n}\left(S_i(\widetilde{\d})S_\blamb(\bfx)\right)
$$
For example
$$
\sigma_3\cap \Omega^\blamb=\pi_{r,n}\left[\left({1\over 6}{\d^3\over \d x_1^3}+{1\over 2}{\d^2\over \d x_1 \d x_2} +{1\over 3}{\d \over \d x_3}\right)S_\blamb(\bfx)\right]
$$

\claim{}  Let $\Delta_\blamb(\exp(t))$ be the Schur determinant as in (\ref{eq:scud}), attached to the exponential formal power series. For example
$$
\Delta_{(3,2,2)}(\exp(t))=\left|\matrix{\displaystyle{1\over 3!}&\displaystyle{1\over 1!}&\displaystyle{1\over 0!}\cr\cr
\displaystyle{1\over 4!}&\displaystyle{1\over 2!}&\displaystyle{1\over 1!}\cr\cr \displaystyle{1\over 5!}&\displaystyle{1\over 3!}&\displaystyle{1\over 2!}}\right|=15
$$
\claim{\bf Corollary.}\label{cor:cor32}
\be
f^\blamb:=|\blamb|!\cdot \Delta_\blamb(\exp(t))
\ee
\proof
Let  $d:=|\blamb|$.  Then
$$
f^\blamb={\d^d S_\blamb(\bfx)\over \d x_1^d}=
\left|\matrix{S_{\lambda_1}(\bfx)&S_{\lambda_2-1}(\bfx)&\cdots&S_{\lambda_r-r+1}(\bfx)\cr\cr
S_{\lambda_1-1}(\bfx)&S_{\lambda_2}(\bfx)&\cdots&S_{\lambda_r-r+2}(\bfx)\cr
\vdots&\vdots&\ddots&\vdots\cr
S_{\lambda_1+r-1}(\bfx)&S_{\lambda_2+r-2}(\bfx)&\cdots&S_{\lambda_r}(\bfx)
}\right|\label{eq:dete}
$$
Now $S_i(\bfx)=\displaystyle{x_1^i\over i!}+g_i$, where
$$
g_i:=g_i(x_1,x_2,\ldots,x_i)
$$
is a polynomal in which $x_1$ occurs with degree strictly smaller than $i$. Therefore the determinant occurring in (\ref{eq:dete} can be written as
\be
\left|\matrix{\displaystyle{x_1^{\lambda_1}\over \lambda_1!}+g_{\lambda_1}&\displaystyle{x_1^{\lambda_2-1}\over (\lambda_2-1)!}+g_{\lambda_2-1}&\cdots&\displaystyle{x_1^{\lambda_r+r-1}	\over (\lambda_r+r-1)!}+g_{\lambda_r+r-1}\cr\cr\cr
\displaystyle{x_1^{\lambda_1+1}\over( \lambda_1+1)!}+g_{\lambda_1+1}&\displaystyle{x_1^{\lambda_2}\over \lambda_2!}+g_{\lambda_2}&\cdots&\displaystyle{x_1^{\lambda_r+r-2}\over (\lambda_r+r-2)!}+g_{\lambda_r+r-2}\cr
\vdots&\vdots&\ddots&\vdots\cr\cr
\displaystyle{x_1^{\lambda_1+r-1}\over( \lambda_1+r-1)!}+g_{\lambda_1+r-1}&\displaystyle{x_1^{\lambda_2+r-2}\over (\lambda_2+r-2)!}+g_{\lambda_2+r-2}&\cdots&\displaystyle{x_1^{\lambda_r}\over \lambda_r!}+g_{\lambda_r}
}\right|\label{eq:hugedet}
\ee
Easy manipulations with determinants show that (\ref{eq:hugedet})  can be written as
$$
x_1^{d}\Delta_\blamb(\exp(e^t))+F(x_1,x_2,\ldots,x_{\lambda_1})
$$
where $F$ is a polynomial in which  $x_1$ occurs in degree smaller than $d$.
Therefore
$$
f^\blamb:={\d^d S_\blamb(\bfx)\over \d x_1^d}={\d^d\over \d x_1^d}\Big(x_1^{d}\Delta_\blamb(\exp(t))+F(x_1,x_2,\ldots,x_{\lambda_1})\Big)=d!\cdot\Delta_\blamb(\exp(t))
$$ 
\qed
\claim{\bf Example.} The degree of the Schubert variety
$
\Omega^{(3,2,1)}(F^\bullet)
$
is
$$
f^{(3,2,1)}=6!\left|\matrix{\displaystyle{1\over 3!}&\displaystyle{1\over 1!}&0\cr\cr
\displaystyle{1\over 4!}&\displaystyle{1\over 2!}&1\cr\cr
\displaystyle{1\over 5!}&\displaystyle{1\over 3!}&1
}\right|=16
$$
which is also the number of standard Young Tableaux of shape $(3,2,1)$.
\claim{\bf Corollary.} \label{cor:34} {\em 
Let $Y(\blamb)$ denote the Young tableau of the partition $\blamb$. Let $h(x)$ denote the {\em hook length} of the box $x$ of  $Y(\blamb)$. Then
$$
\prod_{x\in Y(\blamb)} h(x)={1\over \Delta_\blamb(\exp(t))}
$$
}

\noindent
\proof As in \cite{smirnov}, the degree of a Schubert variety coincides with the number of standard Young tableaux of shape $\blamb$. Then one invoke the celebrated hook length formula proven by \\qed
\claim{\bf Corollary.}\label{cor:35}
\be
\sum_{d\geq 0}{t^d\over d!}\sum_{\blamb\in\Pcal_r\,|\,|\blamb|=d}f^\blamb s_\blamb(z_1,\ldots,z_r)=\exp(t(z_1+\cdots+z_r))
\ee
as a consequence, for each $d\geq 0$, the Pl\"ucker coordinates of the symmetric polynomial $e_1(\bfz_r)^r$ are the degree of the Schubert varieties $\Omega^\blamb(F^\bullet)$
$$
(z_1+\cdots+z_r)^d=\sum_{\blamb}f^\blamb s_\blamb
$$
\proof
Let $\bfz_r;=(z_1,\ldots,z_r)$  and consider the generating function
$$
\sum_{\blamb\in\Pcal_r}X^r(\blamb)s_\blamb(\bfz_r)
$$
of the basis $X^r(\blamb)$ of $B_r$.
By \cite{BeCoGaVi}, this is equal to
$$
\sigma_+(z_1,\ldots,z_r)X^r(0)=\exp\left(\sum_{i\geq 0}x_ip_i(\bfz_r)\right)X^r(0)
$$
Then
\begin{eqnarray*}
\sum_{d\geq 0}{t^d\over d!}\sum_{|\blamb|=d} f_\blamb X^r(0)s_\blamb(\bfz_r)&=& \exp\left(t{\d\over \d x_1}\right)_{|t=0}\exp\left(\sum_{i\geq 0}x_ip_i(\bfz_r)\right)X^r(0)\cr\cr\cr
&=&\exp((x_1+t)p_1(\bfz_r))_{\bfx=0}=\exp(tp_1(\bfz_r))
\end{eqnarray*}
which concludes the proof because $p_1(\bfz_r)=e_1(\bfz_r)$.
\qed

\noindent
More generally, one can consider Schur polynomials in infinitely many indterminates $\bfz:=(z_1,z_2,\ldots)$. We apply  the translation operator along the $x_1$ diiection to the generating function 
$$
S_\blamb(\bfx,\bfz)=\exp(\sum_{i\geq 0}x_ip_i(\bfz))
$$
of the basis elements of $B$, obtaining:
$$
\exp\left(t{\d\over \d x_1}\right)S_\blamb(\bfx,\bfz)=\exp\left(\exp(t+x_1)p_1(z)+\sum_{j\geq 2}x_jp_j(\bfz)\right)
$$
Evaluating at $\bfx=0$ one gets
$$
\sum_{d\geq 0}{t^d\over d!}\sum_{\blamb\vdash d}S_\blamb(\bfx)=\exp(tp_1(\bfz))
$$
\section{Review on Schubert Derivations}\label{sec:sec4}
Consider  the following vector spaces over the rationals: 
\be
\Vcal:=\QQ[X^{-1},X],\qquad V=V_\infty=\QQ[X],\qquad V_n:={V\over X^n\cdot V}\label{eq:spaces}
\ee
alongwith their restricted duals

$$
\Vcal^*=\bigoplus_{i\in\ZZ}\QQ\cdot \d^i,\qquad V^*:=\bigoplus_{i\geq 0}\QQ
\cdot X^i,\qquad \Vcal^*_n=\bigoplus_{0\leq i< n} \QQ\cdot \d^j.
$$
where $\d^j$ stands for the unique linear form on $\Ucal$ such that  $\d^j(X^i)=\delta^{ij}$.
There is  a natural chain of inclusions
$$
V_n\hookrightarrow V\hookrightarrow \Vcal
$$
where the first map is the natural section $X^i+(X^n)\mapsto X^i$ associated to the canonical projection $\Vcal\mapsto \Vcal_n$ and the second is by seeing a polynomial as a Laurent polynomial with no singular part.
\claim{\bf Exterior Algebras.} For all $r\geq 0$ and all $\blamb\in\Pcal_r$ let
\be
X^r(\blamb)=X^{r-1+\lambda_1}\w\cdots\w X^{\lambda_r}
\ee
The exterior algebra of $V_n$ ($n\in\NN\cup\{\infty\}$) is:
$$
\bw V_n=\bigoplus_{r\geq 0}\bw^r\V_n\qquad\mathrm{where}\qquad 
\bw^r\V_n=\bigoplus_{\blamb\in\Pcal_{r,n}}\QQ\cdot X^r(\blamb)\qquad 
$$
For each $m\in\ZZ$, let $\wX^{m+(0)}:=X^{m}\w X^{m-1}\w X^{m-2}\w \cdots$.
The fermionic Fock space of charge $0$ is $F_0:=\bigoplus_{\blamb\in\Pcal}\QQ\cdot \wX^{0+\blamb}$ where  
\be 
\wX^{0+\blamb}:=X^{\lambda_1}\w X^{-1+\lambda_2}\w X^{-2+\lambda_2}\w \cdots\w X^{r-1+\blamb}\w \wX^{r-1+(0)}
\ee
expression which does not depend on $r\geq \ell(\blamb)$.
The {\em boson-fermion correspondence} can be phrased by saying that
$F_0$ is an invertible $B$-module generated by $\wX^0:=\wX^{0+(0)}$, such that
$$
\wX^{0+\blamb}=S_\blamb(\bfx)\wX^0
$$

 Let now $\Ucal$ denote anyone of the spaces listed in (\ref{eq:spaces}).

\bclm{\bf Definition.} {\em A {\em Hasse-Schmidt} (HS) derivation on $\bw \Ucal$ is a $\QQ$-linear map $\Dcal(z):\bw \Ucal\sra \bw \Ucal\llb z\rrb$ such that
\be
\Dcal(z)(u\w v)=\Dcal(z)u\w \Dcal(z)v,\qquad \forall u,v\in\bw \Ucal\label{eq:HSD1}
\ee
}
\eclm
Writing $\Dcal(z)=\sum_{i\geq 0}D_iz^i\in \End_\QQ(\bw \Ucal)\llb z\rrb$, equation (\ref{eq:HSD1}) is equivalent to
\be
D_j(u\w v)=\sum_{i=0}^jD_iu\w D_{j-i}v.
\ee
If  $A\in \End(\Ucal)$,  denote by $\delta(A)\in \End_\QQ(\bw\Ucal)$ the unique derivation of $\bw\Ucal$ such that
\be
\delta(A)u=A\cdot u, \qquad \forall u\in \Ucal=\bw^1\Ucal.
\ee

\bclm{\bf Proposition.}  {\em The plethistic exponential
\be \Dcal^A(z)=\Exp(\delta(A^i)z)=\exp\left(\sum_{i\geq 1}{1\over i}\delta(A)z^i\right):\bw\Ucal\lra \bw\Ucal\llb z\rrb
\ee
is the unique Hasse-Schmidt (HS) derivation on $\bw\Ucal$,
such that
 $\Dcal^A(z)u=\sum_{i\geq 0}(A^iu)z^i$.
 }
\eclm
\proof Based on  the general fact that the exponential of a derivation of an algebra is an algebra homomorphism.\qed

Abusing notation $X$ will also stand for the endomorphism of $\Ucal$ 
given by $u\mapsto Xu$, which is nilpotent if $\Ucal=V_n$ and $n<\infty$.

\bclm{\bf Definition.} {\em The {\em Schubert derivation} $\sigma_+(z):\bw V_n\sra \bw V_n[[z]]$ is 
\be 
\sigma_+(z):=\sum_{i\geq 0}\sigma_iz^i=\Exp(\delta(X)z).
\ee
Its invese is
\be 
\ovsig_+(z):=\sum_{i\geq 0}(-1)^i\ovsig_iz^i=\Exp(-\delta(X)z).
\ee
}
\eclm
They are clearly the unique HS derivation such that $\sigma_+(z)u=\sum_{i\geq 0}X^iu\cdot z^i$ and  $\ovsig_+(z)u=u-Xu$, for all $u\in V_n$.

\claim{} If $u=\sum_{\blamb\in\Pcal_{r,n}}a_\blamb\cdot X^r(\blamb)\in\bw^rV_n$, we denote by $u^c$ the sum $\sum_{\blamb\in\Pcal^{r,n}}a_\blamb X^r(\blamb^c)$.

\bclm{\bf Proposition.} {\em Let $\sigma_-(z):=\sigma_+(z)^*$ be the  $\langle,\rangle$--adjoint of the Schubert derivation $\sigma_+(z)$, i.e.
\be
\langle \sigma_-(z)u,v\rangle=\langle u,\sigma_+(z)v\rangle.
\ee
Then $\sigma_-(z)=\Exp(\delta(X^{-1})z)$, where $X^{-1}$ is the unique endomorphism of $V_n$ mapping $X^j$ to $X^{j-1}$ if $j\geq 1$ and to $0$ otherwise.  
}
\eclm
\proof To show that $\sigma_-(z)$ is a HS-derivation one first identifies $V_n^*$ with $V_n$ through the isomorphism $u\mapsto \langle u,\cdot\rangle$ and then arguues as in \cite[p.~]{pluckercone}. Then one observes that 
$$
\langle \sigma_-jX^i,X^k\rangle=\langle X^i,\sigma_jX^k\rangle=\langle X^i,X^{j+k}\rangle=\delta^{i,j+k}=\delta^{i-j,k}=\langle X^{i-j},X^k\rangle
$$
which  proves that $\ovsig_{-j}X^i=X^{i-j}$. Thus $\ovsig_-(z)=\Exp(\delta(X^{-1})z)$,  because both sides  restrict to the same endomorphism of  $V_n$. \qed

Recall the notation \ref{notpier}. The $\langle,\rangle$--adjoint $\sigma_-(z)$  of $\sigma_+(z)$ will be called Schubert derivation as well.  The reason is due to:
\claim{\bf Theorem.} {\em Schubert derivations $\sigma_\pm(z)$ satisfy Pieri's rule, i.e.
\be
\sigma_iX^r(\blamb)=\sum_{\bmu\in PF_i(\blamb)}X^r(\bmu)\label{eq:mypier}
\ee
and
\be
\sigma_{-i}X^r(\blamb)=\sum_{\bmu\in PF_{-i}(\blamb)}X^r(\bmu)\label{eq:mypier-}
\ee
}
\proof
Formula \ref{eq:mypier} is, up to the notation, \cite[Theorem 2.4]{SCHSD}. To prove (\ref{eq:mypier-}), due to the fact that  $(X^r(\blamb))_{\blamb\in\Pcal_{r,n}}$ is an orthonormal basis of $\bw^rV_n$, one has
\begin{eqnarray*}
\sigma_{-i}X^r(\blamb)&=&\sum_{\bmu\in\Pcal_{r,n}}\langle \sigma_{-i}
X(\blamb), X^r(\bmu)\rangle X^r(\bmu)\cr\cr
&=&\sum_{\bmu\in\Pcal_{r,n}}\langle X^r(\blamb),\sigma_i X^r(\bmu)\rangle 
X^r(\bmu)
\end{eqnarray*}

Using Pieri's formula (\ref{eq:mypier}) one has
$$
\langle X^r(\blamb),\sigma_i X^r(\bmu)\rangle=\langle X^r(\blamb),\sum_{\bnu\in PF_i(\bmu)}
X^r(\bnu)\rangle
$$
and $\bnu\in PF_i(\bmu)$ if and only if
$
|\bnu|=|\bmu|+i\qquad\mathrm{and}\qquad \nu_1\geq \mu_1\geq\cdots\geq\nu_r
\geq \mu_r.
$
Thus
$$
\langle X^r(\blamb),\sigma_i X^r(\bmu)\rangle=\delta_{\blamb,\bnu}
$$
i.e. the only non zero coefficients are those for which $\nu_1=\lambda_1,
\ldots,\nu_r=\lambda_r$, which are then the summands of $\sigma_{-i}X^r(\blamb)$.
\qed

%

\section{Proof of Theorem 1 and one generalization}\label{sec:proof}

By Proposition~\ref{prop:prop24}.  the product $S_i(\bfx)S_\blamb(\bfx)$ obeys  Pieri's formula. With respect to the inner product $\langle,\rangle$ for which $(S_\blamb(\bfx))_{\blamb\in\Pcal}$ is an orthonormal basis of $B$,
one has
\be
\left\langle \exp(\sum_{i\geq 0}x_it^i)S_\blamb(\bfx),S_\bmu(\bfx)\right\rangle=\left\langle S_\blamb(\bfx),\exp\left(\sum_{i\geq 1}{t^i\over i}{\d\over \d x^i}\right)\right\rangle\label{eq:s+s-ad}
\ee
Because of (\ref{eq:s+s-ad}), it follows that the coefficients $S_i(\widetilde{\d})$ of $\exp\left(\displaystyle{\sum_{i\geq 1}{t^i\over i}{\d\over \d x^i}}\right)$ satisfy the dual Pieri formula as in (\ref{eq:mypier-}).
 Therefore
\begin{eqnarray*}
\pi_{r,n}(S_i(\widetilde{\d})S_\blamb(\bfx))&=&\pi_{r,n}\big(\sum_{\bmu\in PF_{-i}(\blamb)}S_\blamb(\bfx)\Big)=\sum_{\bmu\in PF_{-i}(\blamb)}\pi_{r,n}(S_\blamb(\bfx))\cr\cr\cr
&=&\sum_{\bmu\in PF_{-i}(\blamb)}\Omega^\bmu=\sigma_i\cap \Omega^\blamb.
\end{eqnarray*}
\qed

 Let
\be
F(\bfu, \bft):=\left.\exp\left(\sum_{i\geq 1}t_i{\d\over \d x_i}\right)\right|_{x_i=0}\exp\left(\sum_{i\geq 1}u_iS_i(\bfx)\right).
\ee
\claim{\bf Proposition.} {\em 
\be
F(\bfu,\bft)=\exp\left(\sum_{i\geq 1}u_iS_i(\bft)\right),\label{eq:integrals}
\ee
where $S_i(\bft)$ is the Schur polynomial in the variable $\bft$.
}

\smallskip
\proof
It amounts to straightforward manipulation with the Taylor formula.\qed

\claim{\bf Remark.} Formula (\ref{eq:integrals}) is the generating function of the ``integrals'' of product of special Schubert cycles. 
Putting  $h_i:=S_i(\blamb(\bfx))$, it is the generating functions of
$$
\left({\d\over \d x_1}\right)^{|\bmu|}h_\bmu
$$
where if $\bmu=(\mu_1,\ldots,\mu_r)\in\Pcal_r$, one sets
$
h_\bmu:=h_{\mu_1}\cdots h_{\mu_r}
$.

\claim{} A remarkable special case is obtained by setting $t_1=t$, and $t_j=0$ for all $j\geq 2$:
\be
F(\bfu, t)=\exp\left(\sum_{i\geq 1}u_i{t^i\over i!}\right).
\ee
If $\bmu=(1^{m_1}\cdots r^{m_r})$ it is easy to see that
$$
\left({\d\over \d x^1}\right)^{|\bmu|}h_\bmu={(m_1+2m_2+\cdots +rm_r)!\over 1!(2!)^{m_2}\cdots(r!)^{m_r}}
$$
because it is merely the coefficient of $\displaystyle{t^{m_1+\cdots+r m_r}\over (m_1+\cdots+r m_r)!}$
in the expansion of $F(\bfu,t)$. In particular
\be
h_\bmu=\sum_{\blamb\vdash m_1+\cdots+rm_r}<h_\bmu,S_\blamb(\bfx)>S_\blamb(\bfx)\label{eq:itera}
\ee
from which, iterating $(m_1+2m_2+\cdots+rm_r)$ times the derivative with respect to $x_1$ of (\ref{eq:itera}):
\be
{(m_1+2m_2+\cdots +rm_r)!\over 1!(2!)^{m_2}\cdots(r!)^{m_r}}=\sum_{\blamb\vdash m_1+2m_2+\cdots+rm_r}S_\bmu(\widetilde{\d})S_\blamb(\bfx)\cdot f^\blamb\label{eq:lasf}
\ee
which generalizes (\ref{eq:square}). Special cases of (\ref{eq:lasf}) are studied explicitly in \cite{BeGa2}.

\medskip
\noindent
{\bf Acknowledgment.} For supplying the opportunity to give the talk and for logistic and financial support, the Scientific and Organizing Committee of INPANGA 2020[+1] are warmly  acknowledged.  Thanks are due Piotr Pragacz and  Pietro Pirola for helpful  discussions.
Work sponsored by INPAN and, partially, by  Finanziamento Diffuso della Ricerca, no. 53 RBA17GATLET del Politecnico di
Torino.

%
%
%

\bibliographystyle{amsplain}

\begin{thebibliography}{99}

\bibitem{BeThesis} O.~Behzad, {\em Hasse-Schmidt Derivations and Vertex Operators on Exterior Algebras}, Ph.D. Thesis, Institute for Advanced Studies in Basic Sciences, 2021.

\bibitem{BeGa2} O.~Behzad, L.~Gatto, {\em Integrals on Ferminic Fock Spaces}, in progress (2021).

\bibitem{BeCoGaVi}
O.~Behzad, A.~Contiero, L.~Gatto, and R.~Vidal Martins, \emph{Polynomial
  representations of endomorphisms of exterior powers}, Collect. Math. (2021), \href{https://link.springer.com/article/10.1007/s13348-020-00310-5}{https://doi.org/10.1007/s13348-020-00310-5.}

\bibitem{Fulyoung}
W.~Fulton, \emph{Young tableaux}, London Mathematical Society Student Texts,
  vol.~35, Cambridge University Press, Cambridge, 1997, With applications to
  representation theory and geometry. 
\bibitem{Ful}
William Fulton, \emph{Intersection theory}, Springer, 1984.

\bibitem{SCHSD}
L.~Gatto, \emph{Schubert calculus via {H}asse-{S}chmidt derivations}, Asian J.
  Math. \textbf{9} (2005), no.~3, 315--321. 

\bibitem{HSDGA}
L.~Gatto and P.~Salehyan, \emph{Hasse-{S}chmidt derivations on {G}rassmann
  algebras}, IMPA Monographs, vol.~4, Springer, [Cham], 2016, With applications
  to vertex operators.

\bibitem{pluckercone}
L.~Gatto and P.~Salehyan, \emph{On {P}l\"ucker equations characterizing
  {G}rassmann cones}, Schubert varieties, equivariant cohomology and
  characteristic classes\, ---\,{IMPANGA} 15, EMS Ser. Congr. Rep., Eur. Math.
  Soc., Z\"urich, 2018, pp.~97--125,
  \href{https://ga.org/pdf/1603.00510.pdf}{\tt{arXiv:1603.00510}}. 
  
\bibitem{gln}
\bysame, \emph{The cohomology of the {G}rassmannian is a {$gl_n$}-module},
  Comm. Algebra \textbf{48} (2020), no.~1, 274--290. 

\bibitem{SDIWP}
\bysame, \emph{Schubert derivations on the infinite exterior power}, Bull.
  Braz. Math. Soc., New Series \textbf{xxx} (2020), no.~1, 2s. 

\bibitem{gillespie}
M.~Gillespie, \emph{Variations on a theme of Schubert calculus},
  pp.~115--158, Springer International Publishing, Cham, 2019.

\bibitem{thrall}
J.~S.~Frame, G. de B.~Robinson, R.~M.~Thrall  \emph{The hook graphs of the
  symmetric groups}, Canad. J. Math. \textbf{6} (1954), 316--324. 

\bibitem{BeGa}
O.~Behzad, L.Gatto \emph{{Bosonic and Fermionic Representations of
  Endomorphisms of Exterior Algebras}},
  {\href{https://arxiv.org/abs/2009.00479}{\tt ArXiv:2009.00479}}, 2018.

\bibitem{smirnov}
E.~Smirnov, \emph{Grassmannians, flag varieties, and {G}elfand-{Z}etlin
  polytopes}, Recent developments in representation theory, Contemp. Math.,
  vol. 673, Amer. Math. Soc., Providence, RI, 2016, pp.~179--226. 

\end{thebibliography}

\bigskip
\noindent
{\rm Letterio~Gatto}\\
{\tt \href{mailto:letterio.gatto@polito.it}{letterio.gatto@polito.it}}\\
{\it Dipartimento~di~Scienze~Matematiche}\\
{\it Politecnico di Torino}\

\end{document}